\documentclass[10pt]{article}

\usepackage[a4paper, total={6.5in, 9in}]{geometry}
\usepackage{amsmath}
\usepackage{graphicx}
\usepackage{changepage}
\usepackage{mathtools}
\usepackage[english]{babel}
\usepackage[utf8]{inputenc}

\usepackage{fancyhdr} 
\numberwithin{equation}{section}

\usepackage{enumerate}

\newcounter{daggerfootnote}
\newcommand*{\daggerfootnote}[1]{%
    \setcounter{daggerfootnote}{\value{footnote}}%
    \renewcommand*{\thefootnote}{\fnsymbol{footnote}}%
    \footnote[2]{#1}%
    \setcounter{footnote}{\value{daggerfootnote}}%
    \renewcommand*{\thefootnote}{\arabic{footnote}}%
    }

\def \bu {{\hskip -.1in}}

\def \bu {\hskip -.001in}

\def\tttt #1{{\textstyle{#1} }}

\def \Bu {{\bf u}}

\def\tttt #1{{\textstyle{#1} }}

\def \magstep#1 {\ifcase#1 1000\or 1200\or 1440\or 1728\or 2074\or 2488\fi\relax}

\def\la{{\lambda}}

\font\ita=cmssi10  
\font\small=cmr6
\font\title=cmbx10 scaled\magstep2
\font\normal=cmr10 
\font\small=cmr6

\font\bol=cmbx12

\def \-> {\rightarrow}
\def\LL{\big\langle}

\def\RR {\big\rangle}

\def\la {\lambda}

\def \RA {\rightarrow}

\def \sas {\vskip .06truein}
\def\sa{{\vskip .125truein}}

\def\sap{{\vskip .25truein}}

\def\aaa {\alpha}

\def\aa {\alpha}

\def \ses {\,=\,}
\def \sps {\, + \,}

\def \sms {\, - \,}

\def \scs {\, , \,}
\def \ess {\enskip}

\def \ssp {\hskip .25em}
\def \bigsp {\hskip .5truein}
\def \part {\vdash}

\normal

\sap
\def\today{\ifcase\month\or
January\or February\or March\or April\or may\or June\or
July\or August\or September\or October\or November\or
December\fi
\space\number\day, \number\year}

\def \bu {{\hskip -.13in}}

\def \RA {{ \rightarrow }}

\font\small=cmr6
\def \scs {\ssp , \ssp}
\def \ess {\enskip}
\def \ssp {\hskip .25em}
\def \bigsp {\hskip .5truein}
\def \part {\vdash}
\font\title=cmbx10 scaled\magstep2
\font\normal=cmr10 
\def\today{\ifcase\month\or
January\or February\or March\or April\or may\or June\or
July\or August\or September\or October\or November\or
December\fi
\space\number\day, \number\year}

\def \I {{\rm I}}

\def \and {\ess\ess\ess\hbox{and}\ess\ess}
\usepackage{titlesec}

\titleformat{\section}
  {\normalfont\fontsize{12}{15}\bfseries}{\thesection.}{1em}{}
  
\usepackage{mathtools}
\newtagform{noparen}{}{}
\usetagform{noparen}  

\begin{document}

\centerline{\title Fibonacci Polynomials}
\sas 
\centerline{\bol by}
\sas
\centerline{\bol A. M. Garsia and G. Ganzberger \daggerfootnote{Both\ess authors\ess  supported\ess  by\ess  NSF\ess  Grant.}}
\sa

\noindent{\normalfont\fontsize{12}{15}\bfseries Abstract} 
 
 The Fibonacci polynomials $\big\{F_n(x)\big\}_{n\ge 0}$ have been studied  in multiple ways,
[1,7,8,10,14].  In this paper we study them by means of the theory of Heaps of Viennot [12,13]. In this setting our polynomials form a basis $\big\{P_n(x)\big\}_{n\ge 0}$ with 
$P_n(x)$ monic of degree $n$. This given, we are forced to set $P_n(x)=F_{n+1}(x)$. The Heaps setting extends the Flajolet view [4] of the classical theory of orthogonal polynomials given by a three term recursion [3,11]. Thus with Heaps most of the identities for our $P_n(x)'s$ can be derived by combinatorial arguments.
Using the present setting we derive a variety of new identities. We must  mention  that the theory of Heaps is presented here without restrictions. This is much more than needed to deal with the Fibonacci polynomials. We do this to convey a flavor of the power of Heaps.
In the lecture  notes [6] there is a chapter dedicated to Heaps where most of its contents are dedicated to   applications of the theory.

\noindent{\normalfont\fontsize{12}{15}\bfseries Introduction}
 
The  sequence $\big\{P_n(x)\big\}_{n\ge 0}$ is defined by the  recursion
$$
P_{n+1}=x\, P_n(x)+ P_{n-1}(x)
\eqno \I.1
$$
and initial conditions
$$
1)\ess\ess P_{-1}(x)=0,
\bigsp\bigsp
2)\ess\ess P_{0}(x)=1.
\eqno \I.2
$$
Calling $F(x;t)=\sum_{n\ge 0}t^n\,P_n(x)$
the generating function of these polynomials, we see that I.1 and I.2  are equivalent to 
$$
\sum_{n\ge 0} t^{n+1}  P_{n+1}\ses
x t\Big(\sum_{n\ge 0} t^{n}  P_n\Big)\sps
t^2\sum_{n\ge 1} t^{n-1 }   P_{n-1}
\ess\ess\ess\ess \RA \ess\ess\ess\ess
F(x;t) \ses {\frac{1}{1\sms x\, t\,  \sms t^2}}
\eqno \I.3
$$

To give a flavor of the identities we will be able to prove by means of this setting  
we need to define the following sequence of
``{\ita moments}'' (see [15])
$$
\nu_{n}\ses \begin{cases} \frac{(-1)^m}{m+1} \, {2m\choose m} = (-1)^m \, 4^{m+1} \frac{\int_0^1  x^m\sqrt{\frac{1-x}{x}}\, dx}{2 \pi } & \mbox{if } $n=2m$, \\ 
0 & \mbox{otherwise,}
\bigsp \hbox{(for all $n\ge 0$).}
\end{cases}
\eqno \I.4
$$
This given, we will prove that 
$$
P_n(x)= \frac{1}{d_{n-1}}  det \ssp
\begin{pmatrix}
\nu_0 &\nu_1& \cdots &\nu_n \\
\nu_1 &\nu_2& \cdots &\nu_{n+1} \\
\cdots &\cdots  \cdots&\cdots \\
\nu_{n-1} &\nu_n& \cdots &\nu_{2n-1} \\
1 & x &  \cdots &x^n
\end{pmatrix},
d_n= det  
\begin{pmatrix}
\nu_o &\nu_1&\cdots &\nu_n \\
\nu_1 &\nu_2&\cdots &\nu_{n+1} \\
\cdots &\cdots &\cdots&\cdots \\
\nu_{n-1} &\nu_n&\cdots &\nu_{2n-1} \\
\nu_{n} &\nu_{n+1}&\cdots &\nu_{2n} \\
\end{pmatrix}
= (-1)^{\lceil n/2 \rceil}.
\eqno \I.5
$$
We will also prove the identity 
$$
J(x,0,-1)\ses  
1 + \cfrac{1}{1+\cfrac{x^2}{1+\cfrac{x^2}{\cdots}}} = \sum\limits_{m \ge 0} \frac{(-1)^m}{m+1} {2m\choose m} x^{2m}
\eqno \I.6
$$

In the classical theory of polynomial bases 
$\big\{Q_n(x)\big\}$ satisfying a three term recursion the moments are non negative  numbers
(see for instance Favard [3]). This is obviously not always true in I.4. So to prove these identities we need to apply the theory to a suitable substitute. We will eventually need to do so, but first we must give a detailed presentation of the setting of Viennot's ``{\ita monomer-dimer}'' view of continued fractions and derive the tools that we need to prove the validity of our construction of the Fibonacci polynomials. 

We develop here the general classical theory of orthogonal polynomials generated by a three term recursion in the Heap setting of Viennot [12,13]. This is much more than we need to achieve our goal. However, it would be a disservice to the Algebraic Combinatorial audience not to expose them to the general theory of Heaps. A collection of lecture notes  
(see [6]) has recently been submitted for publication by Springer, this  includes
a chapter on Heaps that is completely unrestricted and contains some of the most beautiful applications of the theory.

\vspace{.2cm}

The contents are divided  
into five  sections of various lengths.
\vspace{.4cm} 

\noindent
In section 
 {\bf 1. The ``moment'' scalar product},
 \begin{adjustwidth}{1.3cm}{0cm}
We define the scalar product of two polynomials with real coefficients using ``moments'' and present some results of the classical theory of orthogonal polynomials. The characterization of  these polynomials and their ``three term recursion'' is the high point of this section. 
\end{adjustwidth}
\vspace{.4cm} 

\noindent
In section
{\bf 2. Heaps of monomers and dimers},
 \begin{adjustwidth}{1.3cm}{0cm}
We introduce here the Viennot theory of Heaps of ``monomers'' and ``dimers''
and terminate with a theorem characterizing Motzkin and Dyck paths as pyramids of Heaps of monomers and dimers.
\end{adjustwidth}
\vspace{.4cm} 

\noindent
In section {\bf 3. Moments and Motzkin paths},
 \begin{adjustwidth}{1.3cm}{0cm}
This is the most important section of the paper. We show how the theory of Heaps is related to the Flajolet ``continued fraction'' setting  of the classical theory of orthogonal polynomials defined by a three term recursion.
In Theorems 3.1,  3.2, 3.3 and  3.4 we prove mostly by combinatorial arguments the basic identities of the classical theory that we need to validate our construction of the Fibonacci polynomials. The section ends with two theorems of the classical theory that are of general interest. Their proofs  have been omitted since their results  are not used in the sequel. 
\end{adjustwidth}
\vspace{.4cm} 

\noindent
In section {\bf 4. Heap identities for the Catalan polynomials},
 \begin{adjustwidth}{1.3cm}{0cm}
 We use here a basis $\big\{ Q_n\big\}_{n\ge 0}$ defined by the two  term recursion 
$Q_{n+1}= x\,Q_n -  Q_{n-1}$, we use them as the closest substitute to the Fibonacci polynomials that can be found in the classical theory. This section applies some of the results  of section 3 to obtain identities for the basis $\big\{ Q_n\big\}_{n\ge 0}$.  
\end{adjustwidth}
\vspace{.4cm} 

\noindent
In section {\bf 5. Proofs of Fibonacci polynomials  identities},
 \begin{adjustwidth}{1.3cm}{0cm} 
We use here a basis $\big\{ Q_n(x;\la)\big\}_{n\ge 0}$ as a substitute for the Fibonacci polynomials. These polynomials contain the extra parameter sequence 
$\big\{\la_i\big\}_{i\ge 1}$ and they are
generated by the two  terms recursion
$\,  Q_{n+1}(x;\la)=  xQ_n(x;\la)-\la_n \, Q_{n-1}
(x;\la)\, $. Our first goal here is  
to apply the   results  of section 3 to obtain   the identities satisfied by the  polynomials 
$Q_n(x;\la)$ and all their closely related facts. This done, we transfer by means of the specialization  $\la_i\RA -1$ for all  $i\ge 1$ to obtain  
identities satisfied by the Fibonacci polynomials.
\end{adjustwidth}
\vspace{.4cm}

An important fact that needs to be mentioned here is that, in the classical theory, the moments are given numerically. In particular the parameters $\la_i$ and $c_i$ occurring in the three term recursion are real numbers and the 
$\la_i$ are actually positive. In the Flajolet setting the parameters $\, \la_1,\la_2,\la_3,\ldots\, $ and 
$c_0,c_1,c_2,\ldots$ can be commutative indeterminates. To introduce the initial Viennot  setting we are forced to view 
$\, \la_1,\la_2,\la_3,\ldots\,$ and 
$\, c_0,c_1,c_2,\ldots\, $ as non-commutative indeterminates. When Viennot passes from his original setting to the Flajolet continued fraction setting  the coefficients $\la_i$ and $c_i$
occurring in the recursion  must be allowed to commute.
\vspace{.4cm} 

In this  writing we found that it is more convenient  to  permit this flexibility of point of view with the proviso to make clear the point of view that is being adopted at the very least 
from the context. Basically, if we start with the {\ita moments} then each parameter $\la_n$ and $c_n$ has a formula in terms of the scalar product of elements of the orthogonal basis. The flexibility we adopted permits stating and proving that the moments and the coefficients of the orthogonal basis are polynomials in the  variables $\, \la_1,\la_2,\la_3,\ldots\, $ and 
$\, c_0,c_1,c_2,\ldots\, $. This flexibility allows even the Rogers-Ramanujan continued fraction to be an application of the theory of Heaps. This is one of the examples treated in the chapter on Heaps in [6].
\pagebreak

\sa

\noindent{\bol 1. The ``moment" scalar product}\sas

 In the classical theory (see [3,10]) the scalar product of two polynomials with real coefficients
 \hbox{$A(x)=\sum_{r=0}^{d_a}a_{r} x^r\, $ and
$\, B(x)=\sum_{s=0}^{d_b}b_{s} x^s\,$ }  is defined by setting
\vskip -.1in
$$
\LL A,B \RR_\aaa\ses \int_{-\infty}^{+\infty}\bu A(x)B(x)\,d\aaa\ses \sum_{r=0}^{d_a}\, \sum_{s=0}^{d_b}
\,\, a_{r}\,b_{s}\,\, \mu_{r+s}^\aa
\, ,
\eqno 1.1
$$
where $\aa(x)$ is a weekly increasing function increasing from $0$ to $1$ 
in a finite interval.  The scalar product in 1.1 is well defined since the matrix
$A_n=\|\mu_{i+j}\|_{i,j=1}^n$ has positive eigenvalue for every $n\ge 0$  (see section 5 for a proof). The definition in 1.1 also shows that all we need are the  ``{\ita moments}''  
\vskip -.1in
$$
\mu_n^\aa \ses \int_{-\infty}^{+\infty}\bu
x^n \,d\aaa .
\eqno 1.2
$$

  The existence of a  measure  giving 
$d\aaa$ and $\mu_n^\aa$ by integration is classically referred to as the {\ita moment problem} and it may be of considerable
analytical difficulty depending on what properties the measure  is required to satisfy. 
Our interest here lies on the nature of the 
relations between the sequence of orthogonal polynomials $\{Q_n\}_{n\geq 0}$, with  respect to the scalar product in 1.1, the sequence of moments
$\{\mu_n^\aaa\}_{n\geq 0}$ and two additional sequences  $\{c_n\}_{n\geq 0}$, $\{\la_n\}_{n\geq 1}$. It develops
that these relations may be beautifully expressed by means of the theory of continued fractions. 
What is remarkable about Flajolet's contribution to this subject is to have noticed that many identities of the classical theory can be established by  combinatorial methods.
The corresponding identities are the contents of the following sequence of theorems which combine
classical results of Jacobi, Rogers, Stieltjes and others. We start with a result that shows that
each $\mu_n^\aaa$ may actually be expressed as a polynomial in the  $c\ssp 's$ and the $\la\ssp 's$.

We develop here the classical theory by showing
that a  basis of monic polynomials $\big\{Q_n(x)\big\}_{n\ge 0}$ is orthogonal with respect to the scalar product in 1.1, that is
\vskip -.15in
$$
1)\ess\ess Q_n(x)\ses x^n\sps \sum_{k=0}^{n-1}\ssp a_{n,k}\ssp x^k,
\ess \bigsp\ess 
2)\ess\ess \langle \ssp Q_n\scs Q_m\ssp\rangle_\aa\ses 0\ess\ess when
\ess\ess n\neq m .
\eqno 1.3
$$
if and only if it satisfies the following
``{\ita three-term recursion}'' (see for instance Favard [3])
$$
Q_{n+1}=(x-c_n)Q_n - \la_n\ssp Q_{n-1},
\eqno 1.4
$$
with initial conditions
$$
1)\ess\ess Q_{-1}(x)\ses 0\,, \bigsp
2)\ess\ess Q_{0}(x)\ses 1,
\eqno 1.5
$$
where  
$$
\la_n \ses 
{
\LL Q_{n-1}, xQ_n\RR_\aa
\over
\LL  Q_{n-1},  Q_{n-1}\RR_\aa
}
\ ,\ess\ess\ess\ess\ess\ess\ess\ess
c_n \ses 
{\LL Q_{n}, xQ_n\RR_\aa\over
\LL  Q_{n},  Q_{n}\RR_\aa}\ess .
\eqno 1.6
$$
The reason for this is very simple, the orthogonality of the basis $\big\{Q_n(x)\big\}_{n\ge 0}$ implies that every polynomial $P(x)$ of degree $d$ has an expansion of the form
\vskip -.2in
$$
P(x)\ses \sum_{i=0}^d\ssp 
{ \langle\ssp P\scs Q_i\ssp \rangle_\aaa\over \langle\ssp Q_i\scs Q_i\ssp \rangle_\aaa}
\ssp Q_i(x).
\eqno 1.7
$$
\vskip -.2in

\noindent
Since 
$$
\LL xQ_n,Q_i\RR_\aaa\ses\LL Q_n,xQ_i\RR_\aaa\ses 0\bigsp \big(\hbox{for all $i\le n-2$}\big)
\eqno 1.8
$$
it follows that
$$
xQ_n(x)\ses 
{ \langle\ssp xQ_n\scs Q_{n-1}\ssp \rangle_\aaa\over \langle\ssp Q_{n-1}\scs Q_{n-1}\ssp \rangle_\aaa}
\ssp Q_{n-1}(x)
\sps
{ \langle\ssp xQ_n\scs Q_{n}\ssp \rangle_\aaa\over \langle\ssp Q_{n}\scs Q_n\ssp \rangle_\aaa}
\ssp Q_{n}(x)
\sps
{ \langle\ssp xQ_n\scs Q_{n+1}\ssp \rangle_\aaa\over \langle\ssp Q_{n+1}\scs Q_{n+1}\ssp \rangle_\aaa}
\ssp Q_{n+1}(x).
\eqno 1.9
$$ 
Since from 1.3 we derive that
$\LL xQ_n,Q_{n+1}\RR_\aaa=\LL x^{n+1},Q_{n+1}\RR_\aaa$ we can write
$$
xQ_n(x)\ses 
\la_n Q_{n-1}(x)
\sps
c_n Q_{n}(x)
\sps
 Q_{n+1}(x).
$$
or better
$$
Q_{n+1}(x)=(x-c_n)Q_n(x) - \la_n\ssp Q_{n-1}(x)  .
$$

\noindent
 Since from the definition in  1.1 of the scalar product $\LL\scs\RR_\aaa$ it follows from  1.6 that each $\la_n$ is a positive  real number, 
we can clearly see by comparing  I.1 to
1.4
that our Fibonacci polynomials cannot be automatically  absorbed into the theory of Heaps.  So to prove  identities for Fibonacci polynomials we need first to see what comes out of the classical theory  for polynomials constructed from the recursion 1.4 when all  $c_n=0$ and all $\la_n=1$. That is 
$$
Q_{n+1}= xQ_n \sms Q_{n-1}.
\eqno 1.10
$$
In the rest of this paper we  will omit the dependence on $\aaa$ and simply use  $\mu_n$ to represent the moments.
\sa
\vspace{.2cm}

\noindent{\bol 2. Heaps of monomers and dimers}\sas

The monomer-dimer setting of Viennot  makes the classical theory even more combinatorial. 
 But before making definitions it is best to start with an example.
In the figure below we have an instance of a {\ita monomer-dimer heap}
\vspace{.2cm}
$$
\vcenter{\hbox{\includegraphics[width=3.6in]{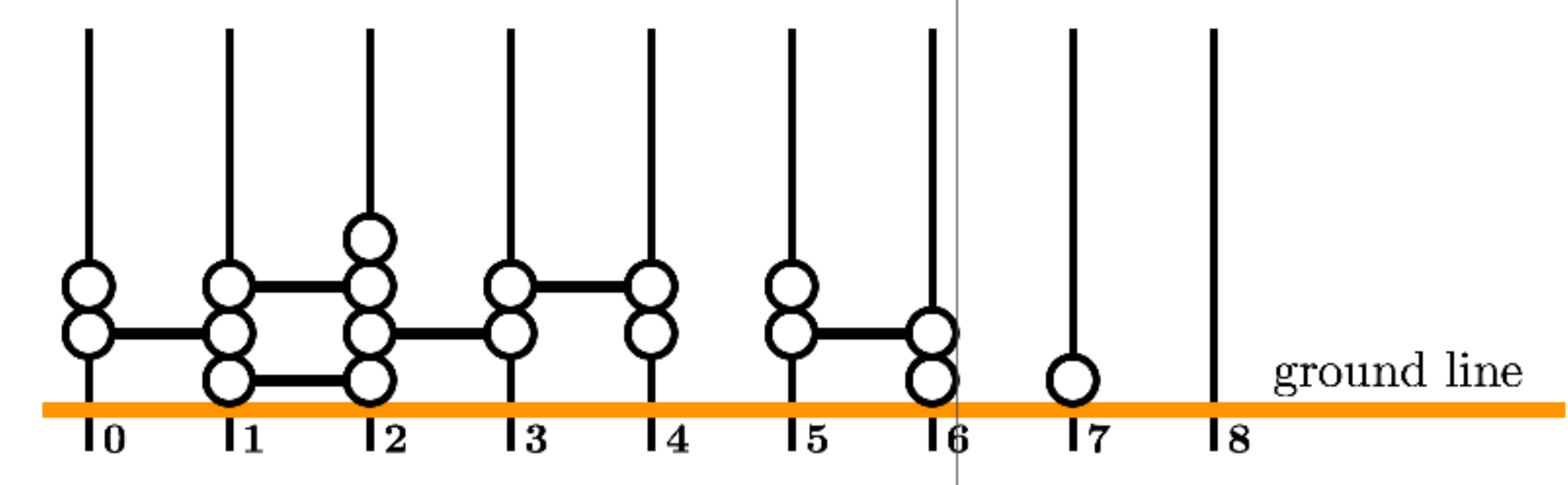}}}
\eqno 2.1
$$
Let us imagine that the vertical lines represent needles and that there are two basic
pieces: 

\sa
\begin{enumerate}
\item A {\ita monomer}, which is a billiard ball pierced along a diameter for  threading by the needles.
\item A {\ita dimer} which consists of two   billiard balls joined by a metal bar. 
\end{enumerate}

\noindent
To put together a heap we simply pick a bunch of monomers and dimers and stack them
on top of each other, threading them by the needles, as indicated in figure 2.1.
As it is depicted there the needles are perpendicular to the ground line at
the points of coordinates $0,1,2,3\ldots$. Heaps of monomers and dimers will be represented
by words in the alphabet
$$
\vspace{.2cm}
{\cal A}\ses \{\ssp m_0,m_1,m_2,\ldots \ssp ;\ssp d_1,d_2,d_3,\ldots \ssp \},
$$
replacing each monomer of ground coordinate $i$ by the letter $m_i$ and each
dimer projecting onto the interval $[i-1,i]$ by the letter $d_i$. The corresponding word is
obtained by processing in this manner the successive pieces of the heap from left
to right within a row, starting from the bottom row and proceeding upwards.
For instance, this procedure applied to the heap of figure 2.1 yields the word
$$
\vspace{.2cm}
w\ses d_2 m_6 m_7 d_1 d_3 m_4 d_6 m_o d_2 d_4 m_5 m_2\ess .
$$
Conversely, given any word $w\in {\cal A}^*$ we can construct a heap by reversing
the process above. That is we read the letters of $w$ from left to right and replace each
$m_i$ by a monomer of ground coordinate $i$ and each $d_i$ by a dimer spanning $[i-1,i]$.
Of course we must also thread the corresponding monomers and dimers down the needles
in the precise succession they are encountered as we read $w$. The final configuration
is obtained by letting the pieces settle as far down as they can. This proceedure
applied to the word $w_1=m_o d_2 m_2 d_1 m_1 d_2 m_3 m_3$ produces the heap
\vspace{.4cm}

$$
\vcenter{\hbox{\includegraphics[width=1.4in]{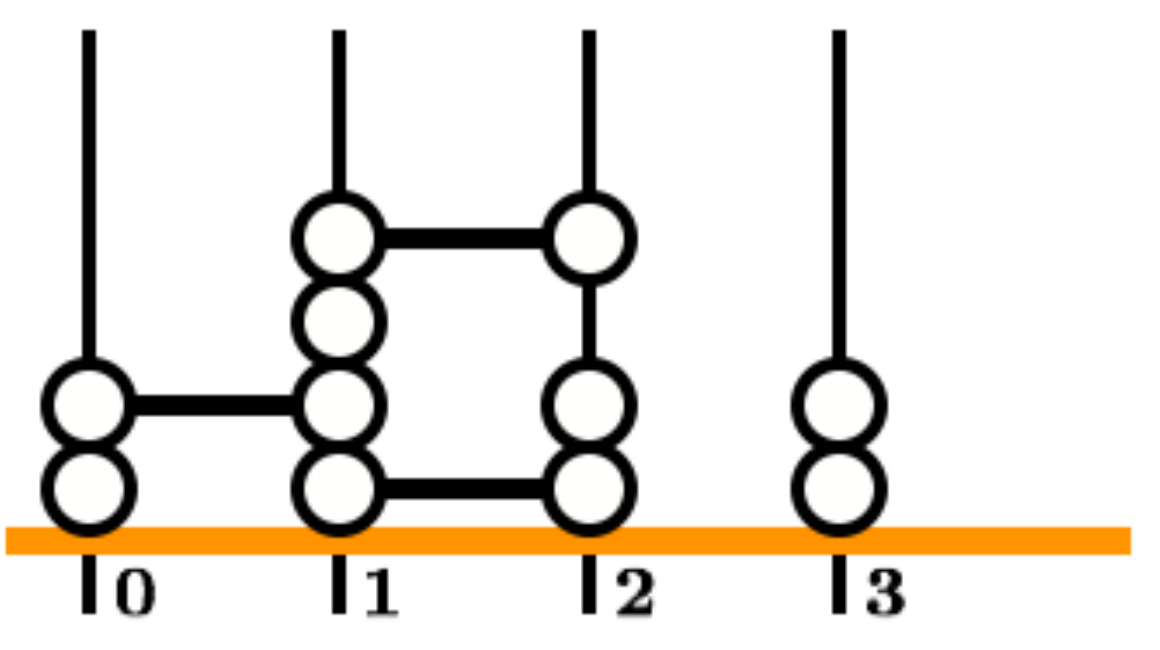}}}
\eqno 2.2
$$

\pagebreak

Now it is easy to see that the word $w_2=m_3 d_2 m_o d_1 m_3 m_1 m_2 d_2 $
produces the same heap. At the same time the word which corresponds to this heap
by the construction given above should be 
$w=m_o d_2 m_3 d_1 m_2 m_3 m_1 d_2$.
Mathematically speaking a ``heap'' should   represent  an equivalence class
of $\cal A$-words. Two words being equivalent if and only if they yield the
same heap. Thus our procedure of constructing the word corresponding to a heap is
just one of the ways of selecting a representative from each equivalence class of
words. 
\sas
\vspace{.2cm}

Imagine now that the needles of figure 2.1 are set into a top and bottom bar as in an 
abacus and we push down the monomer $m_2$. This will result in the configuration
\vspace{.4cm}
$$
\vcenter{\hbox{\includegraphics[width=2.8in]{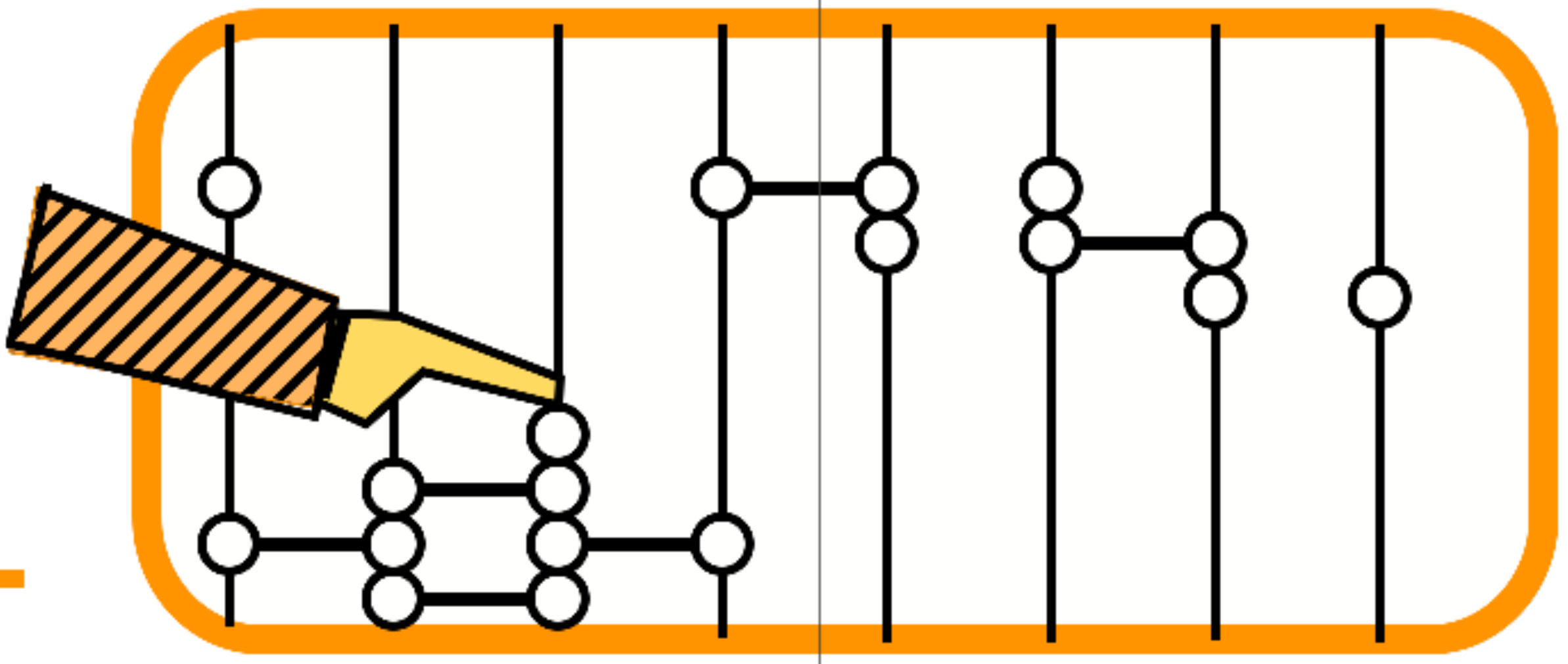}}}
\eqno 2.3
$$
We can see that there are heaps that are brought down by pushing on a single piece.
Such is also the case for the heap obtained by removing the monomers $m_0$ and $m_7$
from the heap in 2.3 and adding on top the dimer $d_5$. Such a heap will be
called a {\ita pyramid} and the top piece the {\ita summit} of the pyramid.

\sa
There is a very simple way of transforming Motzkin paths into heaps
which, as we shall see, has remarkable mathematical consequences. This
transformation is obtained by replacing each {\ita East} step in the path
by a monomer and each {\ita North-East} step by a dimer. We only need one example here to
get across what we have in mind. For instance, carrying out these replacements
(from left to right) on the path on the left yields the configuration in the middle.
The latter is then rotated $90^o$ clockwise so that the pieces settle down to the ground.
This results in the heap given on the right of our display.
\vspace{.4cm}
$$
\vcenter{\hbox{\includegraphics[width=6in]{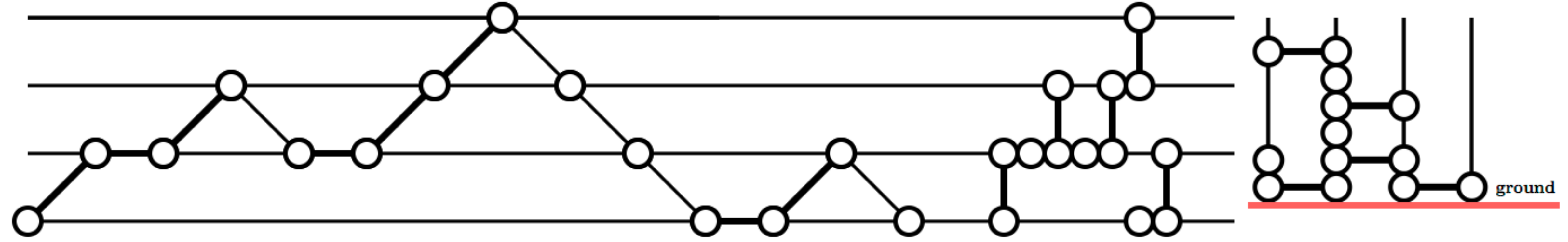}}}
$$
The transformation is even simpler at the word level. In fact, the $a,b,c$-word
corresponding to the Motzkin path is 
$w=a_0c_1a_1b_2c_1a_1a_2b_3b_2b_1c_oa_0b_1$
(here we label the $a,b,c$ letters by the height of the starting point). The sequence of monomers and dimers depicted in the middle of the display is $d_1,m_1,d_2,m_1,d_2,d_3,m_0,d_1$.
 The word of the final heap is 
$w'=d_1d_3m_0d_2m_1d_2m_1d_1$.
We see that to go from the word
$w$ of the path to the word $w'$ of the corresponding heap, we simply read the letters of $w$
from right to left, replace each $a$ by a $d$, each $c$ by an $m$ and remove all 
$b\ssp 's$.
From our example we can easily extract the following
\sas
\vspace{.3cm}

\noindent
{\bol Theorem 2.1}

{\ita Our construction yields a bijection between Motzkin paths and heaps of monomers and
dimers with the following properties.

\begin{enumerate}
\item The image heap is always a pyramid with summit a monomer $m_0$  or a dimer
$d_1$.
\item  Paths whose maximum height does not exceed $n$ correspond to pyramids whose
projection is in the interval $[0,n]$.
\item Dyck paths (no East steps)  are sent into pyramids of dimers with summit $d_1.$
\item If the image pyramid has $d$ dimers and $m$ monomers then the corresponding
path has $2d+m$ steps.
\end{enumerate}}
\pagebreak

\noindent{\bol 3. Moments and Motzkin paths}\sas

Our goal in this section is to derive from the theory  of Heaps all the identities of the classical theory that are  needed to prove the  validity of our construction of the Fibonacci polynomials. We will see 
that using Heaps, the needed classical identities
will  become visually self evident. 
\vspace{.4cm}

To this end we will deal 
here with a more general class of Motzkin paths. These are lattice paths that
proceed by {\ita North-East, East} and {\ita South-East} steps and remain throughout
weakly above the $x$-axis without restrictions on the heights of the starting or ending points.
For instance we give below a Motzkin path that starts at level $5$ and ends at level $2$
\vspace{.5cm}
$$
\vcenter{\hbox{\includegraphics[width=3.8in]{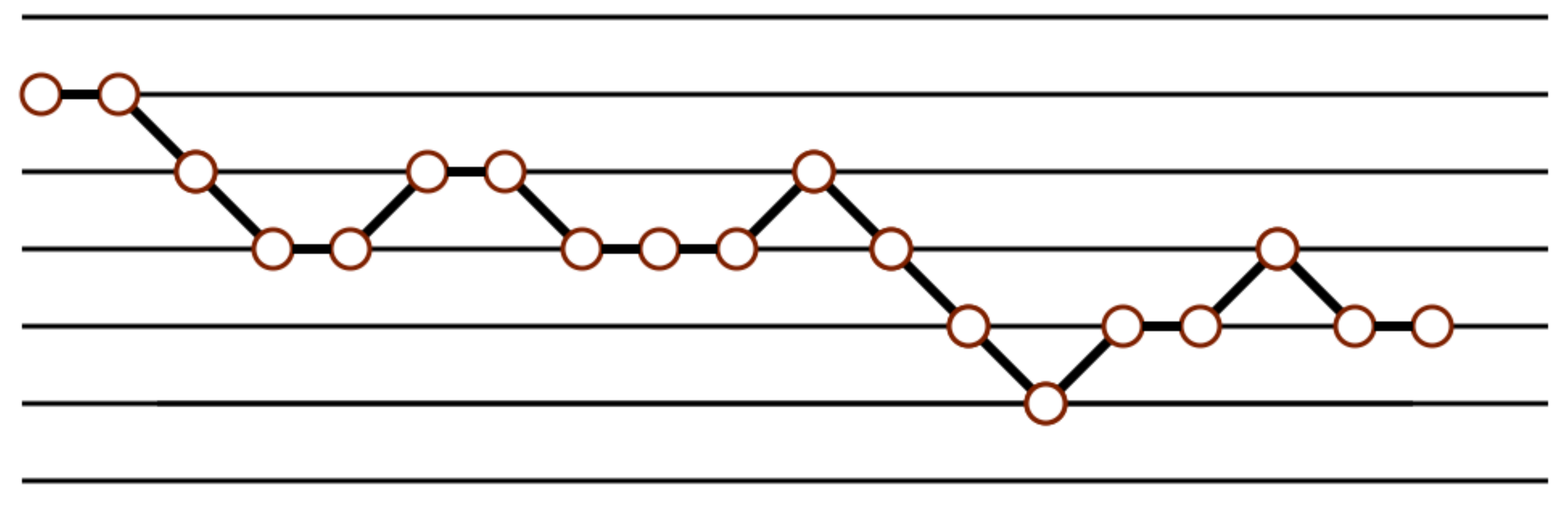}}}
\eqno 3.1
$$
Here and after the collection of Motzkin paths that start at level $r$ and end at level $s$
will be denoted by $\Pi_{r,s}$. To such a path $\pi$ we shall associate a word $w(\pi)$
by replacing (from left to right) each {\ita North-East} edge by an $a_i$,
each  {\ita East} edge by a $c_i$ and each {\ita South-East} edge by a $b_i$,
the subscript $i$ giving the  starting height of the edge. For instance, carrying this
out on the path $\pi$ in 3.1 yields the word
$$
w(\pi)\ses c_5b_5b_4c_3a_3c_4b_4
c_3 c_3 a_3 b_4 b_3 b_2
a_1 c_2 a_2 b_3 c_2. 
$$
Clearly,  we can recover a path from its word. In this context 
it is important to regard $a_i,b_i$ and $c_i$ as sequences of non commuting
variables, sometimes using ``word'' should be the clue. Nevertheless, in other contexts it simplifies our notation to 
allow these letters to commute. We will do so in these contexts as long
as there is no loss. 
\sa
\vspace{.3cm}
\noindent{\bol Proposition 3.1}
\vspace{.4cm}

{\ita Let
$$
{\tilde h}_{n,k}\ses \sum_{\pi\in \Pi_{0,k}^{=n}}\ssp w(\pi)\big|_{a_{i-1}=\la_i\atop b_i=1}, 
\eqno 3.2
$$
where $\Pi_{0,k}^{=n}$ denotes the collection of Motzkin paths that go from height
$0$ to height $k$ in $n$ steps.
Then  the ${\tilde h}_{n,k}$ satisfy the following recursion and initial conditions}
\vspace{.3cm}
$$
a)\ess\ess {\tilde h}_{n,k}\ses 
\la_{k } {\tilde h}_{n-1,k-1}  \sps  c_k {\tilde h}_{n-1,k} \sps    {\tilde h}_{n-1,k+1}
\ess\ess\ess\ess\ess\ess  
b)\ess\ess {\tilde h}_{0,0}=1
\hbox{ with $\, {\tilde h}_{n,k}=0\, $  for $\,n<k$}
\eqno 3.3
$$
\noindent{\bol Proof}

Note that  every path in $\Pi_{0,k}^{=n}$ must come either from height $k-1$ by a {\ita North-East}
step or from height $k$ by an {\ita East} step or from height $k+1$ by a {\ita South-East} step.
This observation yields the recursion
$$
{\tilde h}_{n,k}\ses\la_{k}\, {\tilde h}_{n-1,k-1}  \sps  c_k\, {\tilde h}_{n-1,k} \sps   {\tilde h}_{n-1,k+1}
$$
this proves a) of 3.3.
\vspace{.4cm}

Note further that every path with $n$ edges remains below height $n$ and reaches that 
height only when it consists totally of {\ita  North-East} steps. Thus when $n<k$, 
the right hand side of 3.2 reduces to an empty sum, which forces ${\tilde h}_{n,k}=0$.
Since by definition the word of an empty path reduces to $1$, we also have ${\tilde h}_{0,0}=1$. This proves b) of 3.3 and completes our proof.
\pagebreak

Let now $L$ be 
the formal power series given by the following  summation
$$
L\ses \sum_{\pi\in \Pi_{0,0}}\ssp x^{n(\pi)} \ssp w(\pi),
\eqno 3.4
$$
where $\pi$ is a Motzkin path and $n(\pi)$ denotes the  number of edges of $\pi$.  
\sas
\vspace{.2cm}

 This given,
a moment's reflection should reveal that this formal power series must satisfy
the following identity
$$
L\ses 1 \sps  c_o\, x \ssp L\sps x^2 a_o\ssp (SL)\ssp b_1\ssp L\,,
\eqno 3.5
$$ 
where $S$ denotes the ``shift'' operator that replaces a letter indexed by $i$ 
by the same letter indexed by $i+1$. 
\sas
\vspace{.2cm}

Passing to commutative variables we get 
$$
L\ses
1\, \sps \, L\, \big(\,  c_0x  \sps    x^2  a_0 \, b_1 (SL)\, \big),
\eqno 3.6
$$
equivalently
$$
L\sms L(\,  c_0x  \sps    x^2  a_0 b_1 SL    ) \ses 1.
$$  
This gives
$$
L\big(1\sms     c_0x  \sms    x^2  a_0 b_1 SL   \big) \ses 1
$$
or better
$$
L \ses {1 \over 1\sms     c_0x  \sms    x^2  a_0 b_1 SL}
\eqno 3.7
$$
Successive iterations of 3.7 then yield 
that $L$ must be given by the continued fraction
$$
L=L(x;a,b,c)\ses 
{1\over \displaystyle 1-c_ox-
{\strut a_ob_1 x^2\over \displaystyle 1-c_1x-
{\strut a_1b_2 x^2\over \displaystyle 1-c_2x-
{\strut a_2b_3 x^2\over \displaystyle 1-c_3x-\cdots 
}}}}
\eqno 3.8
$$
From 3.4 it follows that we have 
$$
L\ses\sum_{n\ge 0}x^n  \sum_{\pi\in \Pi_{0,0}(n)} w(\pi),
\eqno 3.9
$$
provided it is understood that  ``$w(\pi)$''    means that it is only a rearrangement of the  non commutative letters of the word of $\pi$.  Here    
   $\Pi_{0,0}(n)$  denotes the collection of paths with $n$ edges from levels $0$ to $0$. 
\sas
\vspace{.2cm}

We can easily  see that any specialization of the sequences $\{a_i\}$,
$\{b_i\}$ that makes $a_{i-1}b_i=\la_i$ reduces   $L(x;a,b,c)$ to   the continued fraction
$$
J(x;c,\la)\ses {1\over \displaystyle 1-c_ox-
{\strut \la_1 x^2\over \displaystyle 1-c_1x-
{\strut \la_2 x^2\over \displaystyle 1-c_2x-
{\strut \la_3 x^2\over \displaystyle 1-c_3x-\cdots 
}}}}
\eqno 3.10
$$
Here and after we shall assume that the equality $L(x;a,b,c)=J(x;c,\la)$ holds true.
\sa

This given,  with the same understanding about the meaning of $``w(\pi)"$,
   we also have
$$
J(x;c,\la)\ses \sum_{n\ge 0} x^n  \sum_{\pi\in \Pi_{0,0}(n)} w(\pi). 
\eqno 3.11
$$
\sas
\pagebreak
Finally we are able to prove the following 
basic facts.
\sas

\vspace{.1cm}

\noindent{\bol Theorem 3.1}

{\ita Let $h_{n,k}=\LL x^n,Q_k\RR$ then 
$$
h_{n,k}=\la_k\, h_{n-1,k-1}+c_k\, h_{n-1,k}+
h_{n-1,k+1},
\eqno 3.12
$$
with the following initial conditions
$$
a)\ess\ess h_{0,0}=1
\bigsp\bigsp
b)\ess\ess h_{n,k}=0 \ess\ess for \ess\ess
n<k\, .
\eqno 3.13
$$
Then
$$
h_{n,k}=J(x;c,\la)\, Q_k^* \Big |_{x^{n+k}},
\eqno 3.14
$$
where
$$
Q_k^*(x)=x^kQ_k(1/x)\, .
\eqno 3.15
$$
In particular, the moment $\mu_n$ of
1.2 is given by the identity}
$$
\mu_n=h_{n,0}
\eqno 3.16
$$
\noindent{\bol Proof}

We begin with the recursion in 1.4 written in the form
$$
xQ_k \ses 
\la_k Q_{k-1} 
\sps
c_k Q_{k} 
\sps
 Q_{k+1}.
\eqno 3.17
$$ 
This given, we derive that
$$
h_{n,k}=\LL x^{n-1},xQ_k\RR\ses
\la_k \LL x^{n-1},Q_{k-1}\RR 
\sps
c_k \LL x^{n-1},Q_{k}\RR 
\sps
\LL x^{n-1}, Q_{k+1}\RR
$$ 
or better
$$
h_{n,k}\ses
\la_k\, h_{n-1,k-1} 
\sps
c_k\, h_{n-1,k} 
\sps
h_{n-1,k+1}, 
$$ 
this proves 3.12.
Since we have $h_{0,0}=\LL x^0,Q_0\RR=1$ because of 1) of 1.3 and we also have  $\LL x^n,Q_k\RR=0$ when $k>n$, we see that the conditions in 3.13 are also satisfied.
We also see that since $Q_0=1$ then 
$h_{n,0}=\LL x^n\RR=\mu_n$, this proves 3.16 as well. The proof of 3.12 and 3.13 shows that 
the sequence $\big\{ {\widetilde h}_{n,k}\big\}_{n\ge k}$ of Proposition 3.1 satisfies the same recursion and the same initial conditions as the  sequence 
$\big\{h _{n,k}\big\}_{n\ge k}$. Thus we  must have the equality
$$
{\widetilde h}_{n,k}\ses h _{n,k}
\bigsp \hbox{(for all $n\ge k$)
}
\eqno 3.18
$$ 
Combining 3.16 with 3.2 for $k=0$ we derive that 
$$
\mu_n\ses \sum_{\pi\in \Pi_{0,0}^{=n}}\ssp w(\pi)\big|_{a_{i-1}=\la_i\atop b_i=1}
$$
and 3.11 becomes
$$
J(x;c,\la)\ses \sum_{n\ge 0} x^n  \mu_n. 
\eqno 3.19
$$
Now since $Q_0=1$ then by 3.13 we also have $Q_0^*=1$. Thus for $k=0$ the identity in 3.14 reduces to 
 $ h_{n,0}=J(x;c,\la)\Big|_{x^n}$, which is true by 3.16 and 3.19. This not only shows   the expansion in 3.19
but also proves the identity
$$
\mu_n\ses \bu \sum_{\pi\in \Pi_{0,0}(n)} w(\pi),\bigsp  \big(\hbox{for all $n\ge 0$}\big),
\eqno 3.20
$$ 
In particular, we derive that each $\mu_n$ is a polynomial in the variables $c_i$ and $\la_i$. For instance the computer gives
$$
\mu_4\ses c_0^4\sps 3\, c_0^2\, \la_1\sps 2\,  c_0\, c_1\, \la_1\sps  c_1^2\, \la_1\sps \la_1^2\sps \la_1\, \la_2.
$$
It remains to verify
3.14 for $ 0<k\le n$. Using 3.19,
1) of 1.3 and 3.15 the identity in 3.14 becomes
$$
h_{n,k}\ses \Big(\sum_{r\ge 0}\mu_r x^r\Big)
\Big(x^n+ \sum_{s=1}^{k-1}a_{k,s}x^{k-s}\Big) \Big |_{x^{n+k}}
\ses \mu_{n+k}\sps \sum_{s=1}^{k-1}a_{k,s}
\, \mu_{n+s}\ses\LL x^n, Q_k\RR.
\eqno 3.21
$$ 
This completes our proof of the theorem.
\sas
\vspace{.2cm}

This shows that the three term recurrence uniquely determines the scalar product with
respect to which the polynomials $Q_k$ are to be orthogonal. 
The essential part of Theorem 3.1 is given by the   equality in 3.14,
which,  in particular, states that the moment sequence $\mu_n$ is given by 3.19.
\sas
According to the original definition 1.1 and 1.2  we have for two polynomials 
$A(x)=\sum_{i=0}^{d_a}a_{i}x^i$ and 
$B(x)=\sum_{j=0}^{d_b}b_{j}x^j$ with real coefficients
$$
\LL x^n \RR\ses \mu_n\,\, \RA
\ess\ess\ess\ess 
\LL A(x),B(x)\RR=
\sum_{i=0}^{d_a}\sum_{j=0}^{d_b}a_{i}b_{j} 
\,\mu_{i+j}.
\eqno 3.22
$$
From 1.3 we have an orthogonal basis of monic polynomials $\big\{Q_n(x)\big\}_{n\ge 0}$ generated by the three term recursion
$$
Q_{n+1}\ses (x-c_n)Q_n-\la_n\, Q_{n-1},
\eqno 3.23
$$
with  initial conditions $Q_{-1}=0$ and $Q_{0}=1$. The orthogonality together with the monic condition implies
$$
Q_n(x)={1\over d_{n-1}} det 
\begin{pmatrix}
\mu_o &\mu_1&\cdots &\mu_n \\
\mu_1 &\mu_2&\cdots &\mu_{n+1} \\
\cdots &\cdots &\cdots&\cdots \\
\mu_{n-1} &\mu_n&\cdots &\mu_{2n-1} \\
1 & x & \cdots &x^n\cr
\end{pmatrix}
\ess\ess\hbox{with}\ess\ess
d_n= det 
\begin{pmatrix}
\mu_o &\mu_1&\cdots &\mu_n \\
\mu_1 &\mu_2&\cdots &\mu_{n+1} \\
\cdots &\cdots &\cdots&\cdots \\
\mu_{n-1} &\mu_n&\cdots &\mu_{2n-1} \\
\mu_{n} &\mu_{n+1}&\cdots &\mu_{2n}
\end{pmatrix}.
\eqno 3.24
$$

This given, we have
\sas

\noindent{\bol Theorem 3.2}
$$
a)\ess\ess {d_n\over d_{n-1}}\ses \la_1\la_2\cdots \la_n
\bigsp ,\bigsp 
b)\ess\ess c_n\ses {\chi_n\over d_n}\sms {\chi_{n-1}\over d_{n-1}}
\eqno 3.25
$$
{\ita with}
$$
a)\ess\ess d_n= det \ssp
\begin{pmatrix}
\mu_o &\mu_1&\cdots &\mu_n \\
\mu_1 &\mu_2&\cdots &\mu_{n+1} \\
\cdots &\cdots&\cdots&\cdots \\
\mu_{n-1} &\mu_n&\cdots &\mu_{2n-1} \\
\mu_{n} &\mu_{n+1}&\cdots &\mu_{2n}
\end{pmatrix},
\ess\ess\ess\ess
b)\ess\ess\chi_n= det \ssp
\begin{pmatrix}
\mu_o &\mu_1&\cdots &\mu_n \\
\mu_1 &\mu_2&\cdots &\mu_{n+1} \\
\cdots &\cdots&\cdots&\cdots \\
\mu_{n-1} &\mu_n&\cdots &\mu_{2n-1} \\
\mu_{n+1} &\mu_{n+2}&\cdots &\mu_{2n+1}
\end{pmatrix}\ess .
\eqno 3.26
$$
\noindent{\bol Proof}

Since there is only one way for a Motzkin path to reach height $n$ from height $0$
in $n$ steps ({\ita North-East}  all the way), formula 3.20 for $k=n$
reduces to (here we use 3.2 and 3.18)
$$ 
\LL Q_n,  Q_n\RR\ses
\LL Q_n,  x^n\RR\ses  
a_oa_1\cdots a_{n-1}\ssp |_{a_{i}=\la_{i+1}}
\ses \la_1\la_2\cdots \la_n\ess .
\eqno 3.27
$$
On the other hand, using 3.24 we get 
$$
\langle Q_n\scs x^n \rangle = {1\over d_{n-1}}\ssp
det\ssp
\begin{pmatrix}
\mu_o &\mu_1&\cdots &\mu_n \\
\mu_1 &\mu_2&\cdots &\mu_{n+1} \\
\cdots &\cdots &\cdots &\cdots \\
\mu_{n-1} &\mu_n& \cdots &\mu_{2n-1} \\
\mu_{n} &\mu_{n+1}& \cdots &\mu_{2n} \\
\end{pmatrix}
\ses 
{d_n\over d_{n-1}}
\eqno 3.28
$$
Combining 3.27 and 3.28 we deduce that
\vspace{.3cm}
$$
{d_n\over d_{n-1}}\ses \la_1\la_2\cdots \la_n,
\eqno 3.29
$$
this proves 3.25 a). To prove 3.25 b) we will proceed purely combinatorially and not 
from the formula in 1.6.
\sas
\vspace{.4cm}
We start by observing that a Motzkin path can reach height $n$ from height $0$ in $n+1$ steps
if and only if it takes $i$ successive {\ita North-East} steps
followed by a single {\ita East} step and then finish up with $n-i$ successive
{\ita North-East} steps (for $i=0,1,\ldots ,n$). By summing all these possibilities (removing non-commutativity) and using the identity in 3.2  (for $k=n$ and $n\RA n+1$) we obtain
\vspace{.4cm}
$$
h_{n+1,n}\ses \sum_{\pi\in \Pi_{0,n}^{=n+1}}\ssp w(\pi)\ssp \Big|_{a_i=\la_{i+1}\atop b_i=1}= 
\big(c_o+c_1+\cdots +c_n\big) \ssp \la_1\la_2\cdots \la_n
\eqno 3.30
$$
Now 3.24 gives
\vspace{.3cm}
$$
\LL x^{n+1},Q_n\RR\ses {\chi_n\over d_{n-1}}
\eqno 3.31
$$
and by combining 3.30 with 3.31 we get
\vspace{.3cm}
$$
\big(c_o+c_1+\cdots +c_n\big) \ssp \la_1\la_2\cdots \la_n
\ses
{\chi_n\over d_{n-1}}.
\eqno 3.32
$$
Finally a use of 3.25 a) yields
\vspace{.3cm}
$$
c_o+c_1+\cdots +c_n 
\ses 
\frac{\chi_n}{d_n},
$$
proving 3.25 b) and completing our proof of Theorem 3.2.
\sas

\vspace{.4cm}

\noindent{\bol Theorem 3.3}

{\ita The matrix $\| h_{n,k}/(\la_1\cdots \la_k)\|$ is the inverse
of the matrix $\|a_{n,k}\|$ of the coefficients of the \hbox{polynomials $Q_n$.}}

\vspace{.4cm}

\noindent{\bol Proof}

Since $\big\{Q_n(x) \big\}_{n\ge 0}$ is an orthogonal basis, we
have
\vspace{.3cm}
$$
\sum_{k=0}^n { \LL y^n, Q_k(y)\RR \over \LL Q_k,Q_k\RR}\, Q_k(x)
\ses x^n.
\eqno 3.33
$$
Expanding the polynomial  $Q_k(x)$ we get  
\vspace{.3cm}
 $$
\sum_{k=0}^n { \LL y^n, Q_k(y)\RR \over \LL Q_k,Q_k\RR}\, 
\sum_{s=0}^k \, a_{k,s}\, x^s
\ses x^n.
$$  
Equating the coefficients of $x^s$ on both sides of this identity we obtain
\vspace{.3cm}
 $$
\sum_{k=0}^n { \LL y^n, Q_k(y)\RR \over \LL Q_k,Q_k\RR}\, 
 a_{k,s}
\ses \chi(s=n).
$$   
In view of the definition of $h_{n,k}$ and 3.27, this identity
is none other than
\vspace{.3cm}
 $$
\sum_{k=0}^n {h_{n,k} \over \la_1\la_2\cdots \la_k}\, 
 a_{k,s}
\ses \chi(s=n), 
\eqno 3.35
$$ 
but that is exactly what the theorem states.
\sas
\pagebreak

This proves the desired result the hard way. The simpler but equivalent way is to  prove instead
\vspace{.2cm}
 $$
\sum_{k=0}^n  a_{n,k}   \, 
 h_{k,s}
\ses \chi(s=n)\la_1\la_2\cdots \la_s. 
$$
This was proved purely combinatorially for $s=n$, (see the proof of 3.27). But for $s<n$, using the definition 
of $h_{k,s}$ given in Theorem 3.1, we have
\vspace{.2cm}
 $$
\sum_{k=0}^n  a_{n,k}   \, 
 h_{k,s}\ses \LL\sum_{k=0}^n  a_{n,k}   \, 
  x^k,Q_s\RR\ses \LL Q_n,Q_s\RR
\ses 0. 
$$

\sa
This completes our derivation of the Heaps identities      
we need  to complete our work on Fibonacci polynomials. However, we will include in this section some additional identities since they might be of interest. Proofs of these identities can be found in the Heaps lecture notes [6].
\sa
\vspace{.2cm}

The finite continued fraction
\vspace{.2cm}
$$
J^{(n)}(x;c,\la)\ses
{1\over \displaystyle 1-c_ox-
{\strut \la_1 x^2\over \displaystyle 1-c_1x-
{\strut \la_2 x^2\over \displaystyle 1-c_2x-
{\strut \cdots \over \displaystyle \cdots 
{\strut \la_n x^2\over \displaystyle 1-c_nx
}}}}}
\eqno 3.36
$$
is usually referred to as the $n^{th}$ {\ita convergent} of $J(x;c,\la)$.
\sas

The next result relates $J^{(n)}(x;c,\la)$ to the  polynomials $Q_n\ssp 's$.
To state it we need further notation.
Given a polynomial $\phi(x;c_o,c_1,\ldots ;\la_1,\la_2,\ldots )$ let $S\phi$ 
denote the polynomial obtained by replacing in $\phi$ each $c_i$ by $c_{i+1}$ and  
each $\la_i$ by $\la_{i+1}$.  
\sa

This given we have:
\sas

\noindent{\bol Theorem 3.4}
{\ita 
\vspace{.2cm}
$$
J^{(n)}(x;c,\la)\ses {SQ_n^*(x)\over Q_{n+1}^*(x)}
\eqno 3.37
$$
where }
$$
1)\ess\ess
Q_n^*(x)=x^nQ_n(1/x),
\bigsp\bigsp 
2)\ess\ess Q_{n+1}^*(x)=x^{n+1}Q_{n+1}(1/x).
\eqno 3.38
$$
\sa 

\noindent{\bol Theorem 3.5}
\vspace{.2cm}
$$
J^{(n)}(x;c,\la)\sms J^{(n-1)}(x;c,\la)\ses {\la_1\la_2\cdots \la_n\ssp x^{2n}\over
 Q_n^*(x) Q_{n+1}^*(x)}
\eqno 3.39
$$

An immediate corollary of 3.39 is that the rational function $J^{(n)}(x;c,\la)$
does converge to $J(x;c,\la)$ at least in the formal power series sense. In fact we
see that the coefficient of $x^n$ in the Taylor series expansion of $J^{(m)}(x;c,\la)$
is the same as that of $J(x;c,\la)$ itself as soon as $2m>n$. This shows that
$\mu_n$ may be directly computed from $J^{({n-1\over 2})}(x;c,\la)$ if $n$ is odd
and from $J^{({n\over2})}(x;c,\la)$ if $n$ is even. In any case, we see that
3.36 defines it to be a polynomial in $c_0,c_1,\ldots ,c_m$ and $\la_1,\la_2,\ldots,\la_m$
where $m$ is the largest integer in $n/2$. Nevertheless, the computation of $\mu_n$
by means of one of the convergents $J^{(m)}(x;c,\la)$, requires (in view of 3.17) the
calculation of the Taylor series of the rational inverse of $Q_{m+1}^*$. To avoid this
step Stieltjes devised an algorithm for the recursive computation of the moments $\mu_n$.
Stieltjes' result is a proof of the identity 3.16.
We proved this identity  
together with other identities in Theorem 3.1.
\sas
 
 We are now ready to apply our results to two classical substitutes for the Fibonacci polynomials.
\pagebreak

\noindent 
{\bol 4. Heap identities for the Catalan polynomials}
\sas 

The simplest substitute is the basis  $\big\{Q_n(x) \big\}_{n\ge 0}$
constructed by the recursion
$$
Q_{n+1}= x\,Q_n -  Q_{n-1},
\eqno 4.1
$$
and initial conditions
$$
1)\ess\ess Q_{-2}(x)\ses 0\,, \bigsp
2)\ess\ess Q_{-1}(x)\ses 1.
\eqno 4.2
$$
Firstly, by comparison with the general case in 2.23 which is 
$$
Q_{n+1}\ses (x-c_n)Q_n-\la_n\, Q_{n-1},
\eqno 4.3
$$
we see that in 4.1 we have
$$
1)\ess\ess  \la_n=1 \ess\ess (\hbox{for all $n\ge 1$}),
\bigsp\hbox{and}\bigsp
2)\ess\ess   c_n=0 \ess\ess  (\hbox{for all $n\ge 0$}).
\eqno 4.4
$$
Combining the definition $h_{n.k}=\LL x_n\scs Q_k\RR$ with
 the result in 3.14 we next obtain
$$
\mu_n\ses \LL x^n\scs Q_0\RR\ses \LL x^n \RR.
\eqno 4.5
$$ 
But  the identity in 3.2  specialized for $k=0$ gives  
$$
\mu_n\ses \sum_{\pi\in \Pi_{0,0}^{=n}}\ssp w(\pi)\ssp \Big|_{a_i=\la_{i+1}=1\atop c_i=0, \, b_i=1}
\eqno 4.6
$$
where $\Pi_{0,0}^{=n}$ denotes the collection of Dyck paths that go from height
$0$ to height $0$ in $n$ steps. 
\sa

Since in  4.6 we are reduced to counting Dyck paths of length $n$ we derive that
$$
\mu_{n}\ses \begin{cases} \tttt{{1 \over m+1}{2\,m \choose m}}  &  \text{if} \hspace{.1cm} $n=2m,$ \vspace{.05cm} \\
0 & \text{otherwise.}
\end{cases}
\eqno 4.7
$$
 For the same reason identity a) of 3.25 forces  $d_n=d_{n-1}$ for all $n\ge 1$. But the second of 3.24 gives  
$d_0=\mu_0=1$. Thus 3.24 yields the following 
 two identities 
 $$
a)\ess\ess Q_n(x)=  \det 
\begin{pmatrix}
\mu_o &\mu_1&\cdots &\mu_n \\
\mu_1 &\mu_2&\cdots &\mu_{n+1} \\
\cdots &\cdots &\cdots&\cdots \\
\mu_{n-1} &\mu_n&\cdots &\mu_{2n-1} \\
1 & x & \cdots &x^n
\end{pmatrix},
\ess\ess\ess\ess\ess\ess\ess\ess 
b)\ess\ess
\det\begin{pmatrix}
\mu_o &\mu_1&\cdots &\mu_n \\
\mu_1 &\mu_2&\cdots &\mu_{n+1} \\
\cdots &\cdots &\cdots&\cdots \\
\mu_{n-1} &\mu_n&\cdots &\mu_{2n-1} \\
\mu_{n} &\mu_{n+1}&\cdots &\mu_{2n}
\end{pmatrix} =1.
\eqno 4.8
$$
Likewise an easy induction, based on b) of 3.25 yields
$$
\det
\begin{pmatrix}
\mu_0 &\mu_1&\cdots &\mu_n \\
\mu_1 &\mu_2&\cdots &\mu_{n+1} \\
\cdots &\cdots&\cdots&\cdots \\
\mu_{n-1} &\mu_n&\cdots &\mu_{2n-1} \\
\mu_{n+1} &\mu_{n+2}&\cdots &\mu_{2n+1}
\end{pmatrix}
\ses \det\begin{pmatrix}1 & 0 \\ 1 & 0 \end{pmatrix}\ses 0.
\eqno 4.9
$$
Finally  we obtain that, in view of 4.7, 3.10 reduces to
$$
L(x;a,b,c)\Big|_{a_{i}b_{i+1}=1\atop c_i=0}
\ses J(x,0,1)\ses   {1\over \displaystyle 1-
{\strut x^2\over \displaystyle 
1-{\strut x^2 \over \displaystyle  
1-{\strut x^2\over \displaystyle \cdots
}}}}\ses \sum_{m\ge 0}x^{2m} \tttt{{1 \over m+1}{2\,m \choose m}}
\eqno 4.10
$$
Due  to this identity  we will call this basis $\big\{Q_n(x)\big\}_{n\ge 0}$
the ``{\ita Catalan Polynomials}''. 

\vspace{.2cm}
 
 \noindent{\bol 5. Proofs of Fibonacci polynomials  identities.} 
\sas
\vspace{.2cm}
The substitute we  use here 
is  the basis $\big\{Q_n(x)\big\}_{n\ge 0}$ 
which satisfies the recursion
$$
Q_{n+1}(x)\ses  x\, Q_n(x)-\la_n\, Q_{n-1}(x),
\eqno 5.1
$$
and initial conditions
$$
1)\ess\ess Q_{-2}(x)\ses 0\,, \bigsp
2)\ess\ess Q_{-1}(x)\ses 1.
\eqno 5.2
$$
Comparing with  the general recursion
$$
Q_{n+1}(x)\ses (x-c_n)Q_n(x)-\la_n\, Q_{n-1}(x),
\eqno 5.3
$$
we see that the only difference is that 
$$
c_n=0,\bigsp\bigsp  \big(\hbox{for all   $n\ge 0$}\big).
\eqno 5.4
$$
Here we start by applying to the  basis $\big\{Q_n(x)\big\}_{n\ge 0}$ all the identities of the general theory. In particular, the corresponding moment sequence will be given by the formula 3.20 for $k=0$ and all $c_n=0$. This is the language of all Dyck paths of length $n$ and we obtain
\vspace{.2cm}
$$
\mu_{n}\ses \begin{cases} 
\displaystyle\sum_{\pi\in \Pi_{0,0}^{=n}}\ssp w(\pi)\ssp \Big|_{a_i=\la_{i+1}\atop b_i=1,  c_i=0}
  &  \text{if } $n=2m$, \vspace{.2cm} \\
0 & \text{otherwise.}\cr
\end{cases}
\eqno 5.5 
$$
To make sure that the meaning of 5.5 is well understood we will illustrate below the case $n=6$
$$
\mu_6\ses 
\la_1\la_2\la_3 \sps
\la_1\la_2\la_2 \sps
\la_1\la_1\la_2 \sps
\la_1\la_2\la_1 \sps
\la_1\la_1\la_1.
\eqno 5.6 
$$

Passing to commutative variables 
from  3.10, 5.4 and 5.5 we get
$$
J(x;0,\la)\ses 
{1\over \displaystyle 1-
{\strut \la_1 x^2\over \displaystyle 
1-{\strut \la_2 x^2 \over \displaystyle  
1-{\strut \la_3 x^2\over \displaystyle \cdots
}}}}\ses \sum_{m\ge 0}\mu_{2m} \,x^{2m}\ess . 
\eqno 5.7
$$
In fact, we can use 5.7 for $\mu_n$  and specialize all the  identities given by 
 3.24,   3.25 and 3.26.  
\sas

Carrying this out yields
\vspace{.2cm}
$$
Q_n(x)={1\over d_{n-1}} det 
\begin{pmatrix}
\mu_o &\mu_1&\cdots &\mu_n \\
\mu_1 &\mu_2&\cdots &\mu_{n+1} \\
\cdots &\cdots &\cdots&\cdots \\
\mu_{n-1} &\mu_n&\cdots &\mu_{2n-1} \\
1 & x & \cdots &x^n
\end{pmatrix}
\ess\ess\hbox{with}\ess\ess
d_n= det 
\begin{pmatrix}
\mu_o &\mu_1&\cdots &\mu_n \\
\mu_1 &\mu_2&\cdots &\mu_{n+1} \\
\cdots &\cdots &\cdots&\cdots \\
\mu_{n-1} &\mu_n&\cdots &\mu_{2n-1} \\
\mu_{n} &\mu_{n+1}&\cdots &\mu_{2n}
\end{pmatrix}.
\eqno 5.8
$$
$$
a)\ess\ess {d_n\over d_{n-1}}\ses \la_1\la_2\cdots \la_n
\bigsp ,\bigsp 
b)\ess\ess  {\chi_n\over d_n}\ses {\chi_{n-1}\over d_{n-1}}
\eqno 5.9
$$
 where we have as in 4.9
 \vspace{.2cm}
$$
\chi_n= det \ssp
\begin{pmatrix}
\mu_o &\mu_1&\cdots &\mu_n \\
\mu_1 &\mu_2&\cdots &\mu_{n+1} \\
\cdots &\cdots&\cdots&\cdots \\
\mu_{n-1} &\mu_n&\cdots &\mu_{2n-1} \\
\mu_{n+1} &\mu_{n+2}&\cdots &\mu_{2n+1}
\end{pmatrix}\ses 0.
\eqno 5.10
$$
\sa 
\pagebreak

Let us now recall that the Fibonacci basis $\big\{P_n(x)\big\}_{n\ge 0}$ satisfies the recursion
$$
P_{n+1}(x)\ses xP_n(x)\sps P_{n-1}(x)
\eqno 5.11
$$
and  initial conditions
$$
1)\ess\ess P_{-1}(x)=0,
\bigsp\bigsp
2)\ess\ess P_{0}(x)=1,
\eqno 5.12
$$
Comparing 5.11, 5.12  with 5.1, 5.2 we see that  to obtain the Fibonacci basis from our present substitute basis is to make the replacements $\la_i\RA -1$ for all $i\ge 1$. Thus all the identities we have established from the classical theory for our substitute basis must remain valid for the Fibonacci basis after this substitution.
\sas
\vspace{.2cm}

Now the first identity is 5.5,  under this substitution (see also 5.6) becomes 
$$
\mu_{n}\ses \begin{cases} 
 \displaystyle\sum_{\pi\in \Pi_{0,0}^{=n}}\ssp w(\pi)\ssp \Big|_{a_i=\la_{i+1}\atop b_i=1, c_i=0}
  &  \text{if }$n=2m$, \vspace{.1cm} \\ 
0 & \text{otherwise,}\cr \end{cases}
\ess\RA\ess\ess\ess 
\nu_{n}\ses \begin{cases} \tttt{{(-1)^m \over m+1}{2\,m \choose m}}  &  \text{if }$n=2m$, \vspace{.4cm} \\
0 & \text{otherwise.}\cr \end{cases}
\eqno 5.13
$$
This proves  I.4. Making the same substitutions on 5.7  gives 
$$
J(x;0,\la)= 
{1\over \displaystyle 1-
{\strut \la_1 x^2\over \displaystyle 
1-{\strut \la_2 x^2 \over \displaystyle  
1-{\strut \la_3 x^2\over \displaystyle \cdots
}}}}=\sum_{m\ge 0}\mu_{2m} \,x^{2m} 
\ess \RA \ess
J(x,0,-1)=  1\sps {1\over \displaystyle 1+
{\strut x^2\over \displaystyle 
1+{\strut x^2 \over \displaystyle  
1+{\strut x^2\over \displaystyle \cdots
}}}}= \sum_{m\ge 0}  \tttt{{(-1)^m \over m+1}{2\,m \choose m}}
\, x^{2m}
$$
This proves  I.6.
\sas
\vspace{.2cm}

Since   
$\mu_0=1$ and 5.9 a)  gives the recursion 
$$
d_n\ses (-1)^n\times d_{n-1},
\eqno 5.14 
$$
which is easily seen to be periodic with period $2$  (after $d_0=1$)
thus  5.9 a) becomes
$$
\vspace{.2cm}
d_n\ses det 
\begin{pmatrix}
\nu_o &\nu_1&\cdots &\nu_n \\
\nu_1 &\nu_2&\cdots &\nu_{n+1} \\
\cdots &\cdots &\cdots&\cdots \\
\nu_{n-1} &\nu_n&\cdots &\nu_{2n-1} \\
\nu_{n} &\nu_{n+1}&\cdots &\nu_{2n}
\end{pmatrix}
\ses (-1)^{\lceil n/2 \rceil},
\eqno 5.15
$$
this proves the second of I.5. From 5.15 and 5.8 we derive that
$$
\vspace{.2cm}
P_n(x)= (-1)^{\lceil {n-1\over 2} \rceil}\, det 
\begin{pmatrix}
\nu_0 &\nu_1&\cdots &\nu_n \\
\nu_1 &\nu_2&\cdots &\nu_{n+1} \\
\cdots &\cdots &\cdots&\cdots \\
\nu_{n-1} &\nu_n&\cdots &\nu_{2n-1} \\
1 & x & \cdots &x^n
\end{pmatrix},
\eqno 5.16
$$
this proves the first of I.5.
\sas

Moreover, using  5.13 the identity in 5.10 becomes 
$$
 det \ssp
\begin{pmatrix}
\nu_0 &\nu_1&\cdots &\nu_n \\
\nu_1 &\nu_2&\cdots &\nu_{n+1} \\
\cdots &\cdots&\cdots&\cdots \\
\nu_{n-1} &\nu_n&\cdots &\nu_{2n-1} \\
\nu_{n+1} &\nu_{n+2}&\cdots &\nu_{2n+1}
\end{pmatrix}\ses 0.
\eqno 5.17
$$

\pagebreak

We terminate with an expansion result in terms of Fibonacci polynomials  which can be stated as a separate 
\sas

\noindent{\bol Proposition 5.1}

{\ita For any polynomial $P(x)$ of degree $d$ we have 
$$
P(x)\ses \sum_{k=0}^d (-1)^{k}\,\langle  P\, ,   P_k  \rangle\,  P_k(x)\, ,   
\eqno 5.18
$$
with a non degenerate scalar product.}

\vspace{.2cm}

\noindent{\bol Proof}

Since $\big\{ P_k(x)\big\}_{k\ge 0}$ is a basis we can certainly have the expansion
$$
P(x)\ses \sum_{k=0}^d\ssp 
{ \langle  P\, , P_k  \rangle \over \langle  P_k\, , P_k  \rangle }
\ssp P_k(x).
\eqno 5.19
$$
However from 3.27 we derive that
$$
\langle  P_k\,, P_k  \rangle\ses (-1)^k,
$$
this proves 5.18. To show that the quadratic form is non-degenerate, the $(n+1)\times (n+1)$ relevant matrix is 
$$
A_n\ses \big\| \nu_{r+s}\big\|_{r,s=0}^{n}\, .
$$
Since addition is commutative this is a symmetric matrix, thus diagonalizable. In particular its determinant gives the product of the eigenvalues.
But we have seen in 5.15 that $\det(A_n)=(-1)^{\lceil n/2 \rceil}$. so none of these eigenvalues can vanish. This completes our proof of the proposition. 
\sas

We purposely programmed on the computer the expansion in 5.18 to obtain $x^n$. What came out is a rather challenging problem. For instance we got
$$
x^7\ses1\sms 14 P_1\sps 14 P_3 \sms 6 P_5\sps P_7
\ess\ess\ess\ess \hbox{and}\ess\ess\ess\ess
x^8\ses 1\sps 14 P_0\sms 28P_2 \sps 20P_4 \sms 7P_6
\sps P_8\, .
\eqno 5.20
$$
We will leave it as a challenge to prove a general formula giving  these computer generated identities.
\sa 

We terminate this section with a list of what is known and what may be new.
\sas

\noindent
Of course it is well known that the Fibonacci polynomials satisfy the recurrence
$$
P_{n+1}(x)\ses xP_n(x)\sps P_{n-1}(x)\, .
\eqno 5.21
$$
The generating function identity 
$$
\sum_{n\ge 0}t^n P_n(x)\ses  {1\over 1-xt-t^2}\, ,
\eqno 5.22
$$
(as we have already seen in the introduction) is an immediate consequence of 5.21.

The formula
$$
P_n(x)\ses {({x+\sqrt {x^2+4}\over 2})^{n+1}-({x-\sqrt {x^2+4}\over 2})^{n+1}\over {x+\sqrt {x^2+4}\over 2}-{x-\sqrt {x^2+4}\over 2}}
\eqno 5.23
$$
follows from 5.22 by the following identities
$$
1-xt-t^2\ses (1-at)(1-bt)\ses 1-(a+b)t+abt^2, 
\eqno 5.24
$$
thus $a+b=x$ and  $ab=-1$. Solving these two identities for $a$ and $b$  yields 
$$
1)\ess\ess a={x+\sqrt {x^2+4}\over 2},
\bigsp\ess\ess
2)\ess\ess b= {x-\sqrt {x^2+4}\over 2}.
\eqno 5.25
$$

On the other hand from 5.22 and 5.24 we also have
$$
{1\over 1-xt-t^2}\ses
\sum_{r\ge 0}t^r\sum_{n+m=r}a^n b^{r-n}\ses
\sum_{r\ge 0}t^r {a^{r+1}-b^{r+1}\over a-b}.
\eqno 5.26
$$
This proves 5.23.
\sas
\vspace{.2cm}

Now from what  Qi and Guo show  in [15] we can derive that 
$$
\vspace{.2cm}
{1\over n+1}{2n \choose n}\ses
4^{n+1}{  {\int_0^1  x^n\sqrt {1-x\over x}\,  dx \over 2\pi}}.
\eqno 5.27
$$
This amazing identity proves that if we set
$$
\vspace{.2cm}
\aaa(u)\ses  \begin{cases} 
\displaystyle \frac{2}{\pi} \int_{0}^{u} \sqrt{\frac{1-x}{x}} dx
  &  \text{if } u \ge 0, \vspace{.2cm} \\
0 & \text{if } u < 0\cr
\end{cases}
\eqno 5.28
$$
then the unit measure $d\aaa(x)$ has the density
$ { 2\over \pi} \, {  {  \sqrt {1-x\over x}   }}$
and 4.7 becomes (using the notation of section 1)
$$
\vspace{.2cm}
\mu_{n}^\aaa\ses \begin{cases} 4^{n+1}{  {\int_0^1  x^n\sqrt {1-x\over x} \,  dx \over 2\pi}}=\displaystyle\int_{-\infty}^{+\infty}\bu
x^n \,d\aaa  &  \text{if } $n=2m$, \vspace{.2cm} \\ 
0 & \text{otherwise}.
\end{cases}
\eqno  5.29
$$
and the right hand side of 5.13 becomes
$$
\vspace{.2cm}
\nu_{n}\ses \begin{cases} (-1)^m 4^{n+1} \frac{\int_{0}^{1} x^n \sqrt{\frac{1-x}{x}}dx}{2\pi}
 &  \text{if } $n=2m$, \vspace{.2cm} \\
0 & \text{otherwise.}
\end{cases}
\eqno 5.30
$$
This should be a new identity. Likewise all (except I.5 [16])  the identities related to the Fibonacci basis $\big\{P_n(x)\big\}_{n\ge 0}$ stated in the introduction, along with 5.17, should be new. The expansion result stated in Proposition 5.1 should also be new as well as the non degeneracy of the scalar product.
\sas

We should give at least an idea of what was known to the classical people  to base their definition
of the scalar product by means of  moments. To simplify our arguments we will work in a very special setting. 
\sas 
 
 \vspace{.1cm}
 
\noindent{\bol Proposition 5.2}
{\ita Suppose that $f(x)\ge 0$ is a continuous function in the interval $[0,1]$  such  that 
$$
\int_0^1f(x)\,dx\ses 1.
\eqno 5.31
$$
Set
$$
\mu_r\ses \int_0^1x^r \, f(x)\,dx,
\bigsp
\hbox{(for all $r\ge 0$).}
\eqno 5.32
$$
Let for each $n\ge 0$ 
$$
A_n\ses \Big\|\mu_{i+j}\Big\|_{i,j=0}^n,
\eqno 5.33
$$
then $A_n$ has only positive eigenvalues..}

\vspace{.2cm}

\noindent{\bol Proof }

 It will be sufficient to carry out our argument in the $3\times 3$ case.
 \pagebreak
 
 Let $\Bu=[u_1,u_2,u_3]^T$ be an eigenvector of $A_3$ with eigenvalue $\la$, then
$$
A_3\left[\begin{matrix} u_1 \\ u_2 \\ u_3 \end{matrix} \right]\ses\left[\begin{matrix}
1 &\mu_1 &\mu_2 \\ 
\mu_1 &\mu_2  &  \mu_3 \\
\mu_2 & \mu_3 & \mu_4
\end{matrix}
\right]
\left[\begin{matrix} u_1\\ u_2 \\ u_3 \end{matrix} \right]
\ses 
\left[\begin{matrix}
u_1+\mu_1u_2+\mu_2u_2 \\
\mu_1 u_1+\mu_2u_2+\mu_3u_2 \\
\mu_2u_1+\mu_3u_2+\mu_4u_2
\end{matrix}
\right]\ses \la \left[\begin{matrix} u_1\\ u_2 \\ u_3 \end{matrix}\right]
\eqno 5.34
$$
thus 
$$
[u_1,u_2,u_3]^T \, A_3\, \left[\begin{matrix} u_1\\ u_2 \\ u_3 \end{matrix} \right]
\ses\la \,(u_1^2+u_2^2+u_3^2)\ses 
\la\,\|\Bu\|^2
\eqno 5.35
$$
and
\begin{align*}
\lambda ||\textbf{u}||^2 = & \int_0^1 \Bigg(x^0u_0(x^0u_0+x^1u_1+ x^2u_2) +  x^1u_1(x^0u_0+ x^1u_1+x^2u_2) + x^2u_2(x^0u_0+x^1u_1+ x^2u_2)
\Bigg) f(x) d \alpha \\
= & \int_0^1 (x^0u_0+x^1u_1+ x^2u_2)^2 f(x) d\alpha >0.
\tag*{5.36}
\end{align*}
This completes our proof.
\sas 
Another result, shown in section 4, is that the Catalan polynomials are close to the Fibonacci polynomials. We find this fact as well as I.4, namely the identity 
$$
\nu_{n}\ses \begin{cases} \tttt{{(-1)^m \over m+1}{2\,m \choose m}} =(-1)^m\,4^{m+1}{  {\int_0^1  x^m\sqrt {1-x\over x}\,  dx \over 2\pi}} &  \text{if } $n=2m$, \vspace{.2cm} \\
0 & \text{otherwise}
\end{cases}
\bigsp \hbox{(for all $n\ge 0$).}
\eqno 5.37
$$
as somewhat  unexpected. 
\sas

The expansion result and the non degeneracy of the scalar product of Fibonacci polynomials suggests that the classical theory can be extended to include arbitrary  real values for the parameters $\{c_i\}_{i\ge 0}$ and $\{\la_i\}_{i\ge 1}$. 

The non-vanishing of the determinants of all the Hankel matrices of the corresponding moments  should 
remain valid even in this extended case.
\sa

\noindent{\bol Acknowledgment}

We want to express our gratitude to Nolan Wallach for discussions and advice during the research that led to this work.


\end{document}